\documentclass[10pt]{amsart}
\usepackage{geometry}                
\usepackage[parfill]{parskip}   
\usepackage{graphicx}
\usepackage{amssymb}
\usepackage{amsmath,amssymb,tikz-cd}
\usepackage{amsthm}
\usepackage{enumitem, longtable, array}

\usepackage{setspace}

\theoremstyle{plain}
\newtheorem{theorem}{Theorem}[section]
\newtheorem{proposition}[theorem]{Proposition}

\newtheorem{conjecture}[theorem]{Conjecture}

\theoremstyle{definition}
\newtheorem{definition}[theorem]{Definition}
\newtheorem{example}[theorem]{Example}

\theoremstyle{remark}
\newtheorem{remark}[theorem]{Remark}

\title{Big Pictures of Motivic and Classical Homotopy Theories}
\author{Ahmad Rouintan} 
\email{amanooman@gmail.com}

\begin{document}
\setstretch{1.1}

\maketitle
\begin{abstract}
	Motivic homotopy theory is meant to play the role of algebraic topology, in particular homotopy theory, in the context of algebraic geometry. As proved by Oliver Röndigs and Paul Arne Østvær, this theory is closely connected to Voevodsky's triangulated category of motives. A connection that is the motivic analogue of the connection between algebraic topology and homological algebra. In this paper, we try to understand the big picture of motivic homotopy theory and its connection to Voevodsky's motives by comparison to the classical counterpart.
\end{abstract}
\tableofcontents

\section*{Introduction}

Motivic homotopy theory is a fast-rising topic in mathematics that Fabien Morel and Vladimir Voevodsky first introduced in the late 90s. Since then, it has faced constant developments, so understanding it in its entirety can be quite overwhelming. Hence, it would be helpful to have a big picture of it guiding you through your studies on this topic. 

This paper is an expository account of motivic homotopy theory and how it connects to Voevodsky's triangulated category of motives. Our approach is to first understand the big picture of classical homotopy theory - connecting algebraic topology to homological algebra - as this would make the task of understanding the big picture of motivic homotopy theory much easier. Then, having this picture in mind, we will study the motivic setting analogously.

There are a few exceptional expository papers available on motivic homotopy theory, each of which covers some parts of the theory. While this seems to defy the need for writing a new one, the unique task of explaining the big picture has not been tackled in any other expository paper fully. The papers \cite{Levine2008} and \cite{Levine2016} by Marc Levine are brilliant expositions, each covering some parts of the big picture while discussing several other results too. The paper \cite{Weibel2004} by Charles A. Weibel is the closest to this paper, but it does not contain discussions about the classical setting. Our aim is to have a more detailed and slightly different paper compared to \cite{Weibel2004} with a chance to compare it to the classical counterpart. Also, in this paper, every theorem, proposition, etc, is referenced properly to guide the readers through their study of motivic homotopy theory (and even classical homotopy theory) with more ease. Having said that, we are only concerned with results that shape the bigger picture of the theories, so by no means this paper is a complete exposition of the topic. 

As it is virtually impossible to cover all the preliminaries for motivic homotopy theory in a paper, familiarity with algebraic topology, homological algebra, category theory, model category theory, and algebraic geometry is assumed.

Here is how this paper is organized:
\begin{enumerate}
\item
In Section \S \ref{C}, we will discuss the big picture of classical homotopy theory. It will include a picture of algebraic topology, a picture of homological algebra, and the connection between the two. 
\item
In Section \S \ref{M1}, we will discuss the first part of the big picture in the motivic setting regarding motivic homotopy theory. This part is analogous to the picture of algebraic topology in the previous section. We will start with a more general discussion about the homotopy theory of simplicial sheaves on a site and then, we apply it to the category of smooth schemes of finite type over a field of characteristic zero with the Nisnevich topology. This gives us the motivic unstable homotopy theory. Lastly, we build the motivic stable homotopy theory by defining the right notion of spectra.
\item
In Section \S \ref{M2}, we will complete our big picture of the motivic homotopy theory by discussing the unstable and stable homotopy theory of motivic spaces with transfers. This part is analogous to the picture of homological algebra in the first section. As we go through this section, we build a link with the previous section which eventually results in a connection between the motivic stable homotopy category and Voevodsky's triangulated category of motives.
\item
In Section \S \ref{T}, we will briefly discuss the conjectural motivic $t\text{-structure}$, which builds a way from Voevodsky's motives to the conjectural category of mixed motives.
\end{enumerate}

\section{The picture in classical homotopy theory}\label{C}

Algebraic topology and homological algebra - the two main theories when studying (co)homology theories - are closely connected together. In a way, the latter is just the \lq\lq linearization" or \lq\lq abelianization" of the former. This connection is the topic of this section. We will study it and summarize it in a diagram consisting of some categories and functors between them. To do so, we will start by stating the axiomatic definition of reduced cohomology theories and use its axioms to describe important categories involved in algebraic topology. Although this backward approach is not how the theory was built, it serves our goal of drawing a big picture perfectly. After that, we build a bridge to homological algebra using the Dold-Kan correspondence and take analogous steps in this theory. By the time we complete these steps, we will have a big picture of classical homotopy theory, which we will use in order to keep track of the development of motivic homotopy theory in later sections. 

For our discussions to be precise, we must first choose a convenient category of spaces first. This is because the usual category of topological spaces allows pathology that prevents a clean development of homotopy theory. We will not discuss the reasons behind that and refer interested readers to Chapter 5 of \cite{May1999}. Anyway, two choices for the category of spaces that lead to equivalent homotopy theories are
\begin{enumerate}
\item
the category of compactly generated Hausdorff topological spaces, denoted by $\mathbf{Top}$, and
\item
the category of simplicial sets, denoted by $s\mathbf{Set}$.
\end{enumerate}
We will use either of these two categories of spaces after we have discussed the equivalence between the resulting homotopy theories.

Lastly, I should mention that to write this section, I got help from the videos of Peter Arndt's mini-course on abstract and motivic homotopy theory\footnote{available at: https://www.youtube.com/playlist?list=PLYShsxzWr0U3EZobFoA7sJcpjip8zXzui} at the University of Verona as well as the papers \cite{Levine2008} and \cite{Levine2016}. 

\subsection{Reduced cohomology theories}

\begin{definition}\label{rct}
	A \emph{reduced cohomology theory} is a family of functors 
	$$\{E^n: \mathbf{Top}_\bullet^\text{op} \to \mathbf{Ab}\}_{n=0}^{\infty},$$
	from non-degenerately pointed topological spaces to abelian groups, satisfying the following axioms:
	\begin{enumerate}
		\item 
		For every weak equivalence $f: X\to Y$ of topological spaces the induced map
		$$f^n: E^n(Y)\to E^n(X)$$
		is an isomorphism for every integer $n\geq 0$. Recall that a weak equivalence is a map that induces isomorphisms on all homotopy groups.
		\item 
		For every topological space $X$ and every integer $n\geq 0$, there is a natural isomorphism 
		$$\Sigma: E^n(X)\to E^{n+1}(\Sigma X)$$ 
		i.e. reduced cohomology theories are stable with respect to the suspension functor. Recall that the suspension of a space $X$, denoted by $\Sigma X$, is its smash product with $S^1$.
		\item 
		It turns the wedge of spaces into the product of their cohomology groups i.e. 
		$$E^* (\bigvee_{i\in I} X_i) \cong \prod_{i\in I} E^* (X_i).$$
		\item 
		For every cofibration $A\to X$, it induces the following exact sequence:
		$$E^*(X/A)\to E^*(X)\to E^*(A).$$
	\end{enumerate}
\end{definition}
Having this definition in mind, we will analyze its axioms one by one and use them as hints to draw our big picture.

\subsection{Homotopy category of spaces}

The first axiom in Definition \ref{rct} states that reduced cohomology theories turn weak equivalences into isomorphism. Hence, if we turn weak equivalences into isomorphism and do it in a universal way, this axiom becomes trivial. The process of turning a set of morphisms into isomorphisms universally is called localization, which is not always possible by the way. To address this, we use the model category theory which allows us to perform localizations - see \cite{Hirschhorn2002} for a good reference on model category theory. 

\begin{theorem}[\cite{Hirschhorn2002}, Theorem 7.10.10 and 7.10.11]\label{topmodel}
The category $\mathbf{Top}$ with the class of usual weak equivalences and the class of Serre fibrations becomes a model category. Also, the category $\mathbf{Top}_\bullet$ of pointed spaces is a model category with its weak equivalences and fibrations being the usual weak equivalences and Serre fibrations of underlying spaces.
\end{theorem}

As a result of Theorem \ref{topmodel}, we can localize with respect to weak equivalences and get the homotopy category of spaces $\mathcal{H}$ and its pointed version $\mathcal{H}_\bullet$. 

\begin{remark}\label{simset}
	As we mentioned earlier, another choice for the category of spaces is the category of simplicial sets. This category and its pointed version enjoy model structures given by the usual weak equivalences and Kan fibrations - see Theorem 7.10.12 and 7.10.13 in \cite{Hirschhorn2002}. Hence, we get the homotopy category of simplicial sets and its pointed version, and not only these homotopy categories are equivalent to $\mathcal{H}$ and $\mathcal{H}_\bullet$ respectively, but the model categories are Quillen equivalent, too. This is given by the geometric realization functor in one way and the singular complex functor in the other. Therefore, we can substitute topological spaces with simplicial sets whenever needed.
\end{remark}

\subsection{Stable homotopy category of spaces}

For cohomology theories to be stable with respect to suspension we need to turn the suspension functor into an auto-equivalence. One way to do that is to use the Spanier-Whitehead suspension category $\mathrm{SW}$. The idea, naively, is to consider the category of pointed spaces and define
$$\text{Mor}_{\text{SW}}(X, Y) := \text{colim}_q [\Sigma^q X, \Sigma^q Y],$$ 
where $X$ and $Y$ are pointed spaces and by $[-, -]$ we mean the homotopy classes of continuous maps - see \cite{SpanierWhitehead1953} for a more detailed and general discussion. Although this category is triangulated with suspension playing the role of its shift functor, it lacks some properties needed. For example, it does not have infinite coproducts, and hence, does not admit a model structure. Also, the canonical functor from $\mathcal{H}_\bullet$ to this category does not preserve already existing coproducts. The reason why we care so much about infinite coproducts is that not having them means that we can not use the Brown representability theorem, a theorem that plays an important role later on. To address such problems we introduce the category of spectra.

\begin{definition}
	A \textit{spectrum} is a sequence of pointed spaces $E = (E_0, E_1, ...)$ together with a bonding map $\epsilon_n : \Sigma E_n\to  E_{n+1}$ for each integer $n\geq 0$. When the bonding maps are weak equivalences, we call our spectrum an $\Omega\textit{-spectrum}$. A map of spectra $f: E\to F$ is a family of maps $\{f_n: E_n \to F_n\}_{n\geq 0}$, that respect the bonding maps i.e. the diagram 
\begin{equation*}
	\begin{tikzcd}
		\Sigma E_n
		\arrow[r, "{\Sigma f_n}"] 
		\arrow[d, "{\epsilon_n}"] & 
		\Sigma F_n
		\arrow[d, "{{\epsilon}^{'}_{n}}"] \\
		E_{n+1}
		\arrow[r, "F_{n+1}"] &
		F_{n+1}
	\end{tikzcd}
\end{equation*}
commutes for each integer $n\geq 0$. We denote the category of spectra by $\mathbf{Spt}$.
\end{definition}
From a pointed space $X$, we can construct an $\Omega\text{-spectrum}$ $\Sigma^\infty X = (X, \Sigma X, \Sigma^2 X, ...)$ with identity bonding maps. This defines a covariant functor $\Sigma^\infty$ from the category of pointed spaces to the category of spectra. The right adjoint of this functor that sends every spectrum to its zeroth term and maps of spectra to the maps between their zeroth terms is denoted by $\Omega^\infty$. 

In the next step, we should turn the category of spectra into a model category and construct its homotopy category. So, we should first define weak equivalences of spectra. In order to do so, let’s go back to a famous question of stable homotopy groups of spheres and recall the Freudenthal suspension theorem. 

\begin{theorem}[Freudenthal suspension theorem, \cite{Freudenthal1938}]\label{Freu}
For integers $n< k+1$, the suspension homomorphism
\begin{equation*}
		\begin{tikzcd}
			\pi_{n+k}(S^k) 
			\arrow[r, "{\Sigma}"] &
			 \pi_{n+k+1}(\Sigma S^k) = \pi_{n+k+1}(S^{k+1})
		\end{tikzcd}
	\end{equation*}
is an isomorphism.
\end{theorem}

Based on Theorem \ref{Freu}, we can see that the homotopy groups of spheres become stable under suspension when the dimension is large enough. So, we can define the $n\text{-th}$ stable homotopy group of spheres as the colimit 
$$\pi_n^s := \mathrm{colim}_k \pi_{n+k}(S^k).$$ 

We can define stable homotopy groups of spectra in a similar fashion.
\begin{definition}
	The $n\textit{-th}$ \textit{stable homotopy group of a spectrum} $E$ is defined to be the colimit 
	$$\pi_n(E) := \text{colim}_k \pi_{n+k}(E_k),$$ 
	where the maps are given by
	\begin{equation*}
		\begin{tikzcd}
			\pi_{n+k}(E_k) 
			\arrow[r, "{\Sigma}"] &
			\pi_{n+k+1} (\Sigma E_{k}) 
			\arrow[r, "{\epsilon_n}_*"] &
			\pi_{n+k+1}(E_{k+1}).
		\end{tikzcd}
	\end{equation*}
\end{definition}

\begin{example}
	The $n\text{-th}$ stable homotopy group of spheres is just the $n\text{-th}$ stable homotopy group of the spectrum $\mathbb{S} = (S^0, S^1, S^2, ...)$ with the obvious bonding maps. More generally, for a pointed space $X$, its $n\text{-th}$ stable homotopy group is 
	$$\pi_n^s(X) = \pi_n(\Sigma^{\infty} X) = \text{colim}_k \pi_{n+k}(\Sigma^k X) = \text{colim}_k \pi_n (\Omega^k(\Sigma^k X)),$$
	where $\Omega$ is the loop-space functor.
\end{example}

Having defined stable homotopy groups of spectra, we define stable weak equivalences similar to the unstable case.

\begin{definition}
	A \textit{stable weak equivalence} is a map of spectra that induces isomorphisms on all stable homotopy groups.
\end{definition}

There is a nice model structure for spectra, where weak equivalences are stable weak equivalences and cofibrations are those maps of spectra $f: E\to F$ where each $f_n$ is a cofibration in the unstable sense. Therefore, we can localize this category with respect to the stable weak equivalences and get the stable homotopy category, denoted by $\mathcal{SH}$. Recall that the objects of this homotopy category are bifibrant spectra that are just $\Omega\text{-spectra}$. We stop our discussion of spectra here, but, there are different categories of spectra admitting different model structures, each of which has its own advantages. For a full discussion of such topics see \cite{Barns2020}. 

The pair of adjoint functors $(\Sigma^\infty \dashv \Omega^\infty)$ raises to a Quillen adjunction between respective model categories and by taking derived functors to a pair of adjoint functors between $\mathcal{H}_\bullet$ and $\mathcal{SH}$. Using $\Sigma^\infty$, we embed our pointed spaces into $\mathbf{Spt}$ where suspension, which now has become a shift functor, is an auto-equivalence with the obvious inverse and this turns $\mathcal{SH}$ into a triangulated category. 

After passing to the stable homotopy category both of the first two axioms would be trivial, and our picture of classical homotopy theory looks like the following: 

\begin{equation}
	\begin{tikzcd}
		\mathbf{Spc}_\bullet 
		\arrow[r] & 
		\mathcal{H}_\bullet
		\arrow[r, shift left=0.5ex, "{\Sigma^\infty}"] & 
		\mathcal{SH}. 
		\arrow[l, shift left=0.5ex, "{\Omega^\infty}"] 
	\end{tikzcd}
\end{equation}

\subsection{Reduced (co)homology theories as spectra}

\begin{definition}\label{sptcoh}
	Let $E$ and $F$ be two spectra. The $n\text{-th}$ $E\textit{-homology}$ of $F$ is defined to be
	$$E_n(F) := \text{Hom}_{\mathcal{SH}}(\Sigma^n \mathbb{S}, E\wedge F),$$
	and the $n\text{-th}$ $E\textit{-cohomology}$ of $F$ is defined to be 
	$$E^n(F) := \text{Hom}_{\mathcal{SH}}(F, \Sigma^n E).$$
	For a pointed space $X$, its $E\text{-homology}$ and $E\text{-cohomology}$ are respectively $E_*(\Sigma^\infty X)$ and $E^*(\Sigma^\infty X)$.
\end{definition}

These (co)homology theories satisfy the first two axioms in definition \ref{rct}. In \cite{Brown1962}, Brown first proves that such (co)homology theories satisfy the last two axioms and then, introduces his famous representability theorem, which states that every reduced cohomology theory is given by an $\Omega\text{-spectrum}$ by the means of Definition \ref{sptcoh}. Then, in \cite{Whitehead1962}, Whitehead uses the Brown representability theorem to deduce that every reduced homology theory is given by Definition \ref{sptcoh} too. Having defined a homology theory, one needs to convert it to a cohomology theory for finite spectra and then represent that cohomology theory. Therefore, we can informally say that $\mathcal{SH}$ is where reduced (co)homology theories live!

\begin{example}[Singular (Co)homology]
	Let $K(\mathbb{Z}, n)$ be an Eilenberg-MacLane space i.e. the space with the only non-zero homotopy group $\pi_n(K(\mathbb{Z}, n)) = \mathbb{Z}$. The corresponding Eilenberg-MacLane spectrum is defined to be the spectrum 
	$$H\mathbb{Z} = (K(\mathbb{Z}, 0), K(\mathbb{Z}, 1), K(\mathbb{Z}, 2), ...)$$ 
	with its bonding maps given by the adjoint of the weak equivalences $K(\mathbb{Z}, n)\to \Omega K(\mathbb{Z}, n+1)$. This spectrum represents singular homology and singular cohomology with coefficients in $\mathbb{Z}$ i.e. 
\begin{align*}
H^n(X; \mathbb{Z}) \cong \text{Hom}_{\mathcal{SH}}(\Sigma^{\infty} X, \Sigma^n H\mathbb{Z}),\\
H_n(X; \mathbb{Z}) \cong \text{Hom}_{\mathcal{SH}}(\Sigma^n \mathbb{S}, H\mathbb{Z} \wedge \Sigma^{\infty} X).\qedhere
\end{align*}
\end{example}


\subsection{The link between spaces and chain complexes of abelian groups}

In this subsection, we are looking to take similar steps in the context of homological algebra and connect them to our previous discussions. Linking spaces to the chain complexes of abelian groups plays a big part in our picture of classical homotopy theory.
\begin{definition}
	A \textit{chain complex of abelian groups} is a sequence $A = \{A_n\}_{n\in \mathbb{Z}}$ of abelian groups with a family of maps $\{d_n:A_{n+1}\to A_n\}_{n\in \mathbb{Z}}$ that $d^2 = 0$. A connective chain complex of abelian groups is one with non-negative support. 
	A map of chain complexes $f:A\to B$, is a family of maps $\{f_n: A_n\to B_n\}_{n\in \mathbb{Z}}$ for which the diagrams
	\begin{equation*}
		\begin{tikzcd}
			A_{n+1}
			\arrow[r, "d_{n+1}"] 
			\arrow[d, "f_{n+1}"] &
			A_n
			\arrow[d, "f_n"] \\
			B_{n+1} 
			\arrow[r, "d^{'}_n"] &
			B_n
		\end{tikzcd}
	\end{equation*}
	are commutative for every integer $n$. We denote the category of connective chain complex of abelian groups by $C_+(\mathbf{Ab})$.
\end{definition}

\begin{definition}
	Let $A$ be a chain complex of abelian groups. The $n\textit{-th}$ \textit{homology group} of $A$ is 
	$$H_n(A) := \mathrm{ker}(d_n)/ \mathrm{im}(d_{n+1}).$$
\end{definition}

\begin{definition}
	A \textit{quasi-isomorphism} is a map of chain complexes that induces isomorphisms on all homology groups.
\end{definition}

Let $A$ be a simplicial abelian group. We can correspond a connective chain complex of abelian groups $N(A)$ to it by defining 
$$N_m(A) := \bigcap_{i=0}^{m-1} \{\text{ker } \partial_i: A_m\to A_{m-1}\}$$
with its differential given by the remaining face map $d_m:= \partial_m$.

\begin{theorem}[Dold-Kan correspondence, \cite{Weibel1994}, Theorem 8.4.1]
The normalized chain complex functor is an equivalence of categories between the category of simplicial abelian groups $s\mathbf{Ab}$ and the category of connective chain complexes of abelian groups $C_+(\mathbf{Ab})$. Moreover, under this equivalence, simplicial homotopy turns into homology and weak equivalences turn into quasi-isomorphisms. 
\end{theorem}

Using the Dold-Kan correspondence, we will induce a model structure on $C_+(\mathbf{Ab})$ with quasi-isomorphisms playing the role of weak equivalences, but first, let's introduce a model structure on $s\mathbf{Ab}$. Based on Theorem 2.4. in Chapter II of \cite{Quillen1967}, the category $s\mathbf{Ab}$ admits a model structure where its weak equivalences (resp. fibrations) are weak equivalences (resp. fibrations) of underlying simplicial sets. We impose this model structure on $C_+(\mathbf{Ab})$ and denote the resulting homotopy category by $\mathcal{H}_+(\mathbf{Ab})$. Because the forgetful functor $U: s\mathbf{Ab} \to s\mathbf{Set}$ preserves both weak equivalences and fibrations, the free-forgetful adjunction $f_\mathbb{Z} \dashv U$ between $s\mathbf{Set}$ and $s\mathbf{Ab}$ becomes a Quillen adjunction and by taking derived functors, we get an adjoint pair 
\begin{equation*}
\begin{tikzcd}
f_\mathbb{Z}: \mathcal{H}
\arrow[r, shift left=0.5ex] &
\mathcal{H}_+(\mathbf{Ab}): U.
\arrow[l, shift left=0.5ex]
\end{tikzcd}
\end{equation*} 
Moreover, this adjunction factors through the homotopy category of pointed spaces as the free-forgetful adjunction $f_\mathbb{Z} \dashv U$ factors through the free-forgetful adjunction between spaces and pointed spaces. Up to this point, out picture of classical homotopy theory looks like the following:
\begin{equation}
	\begin{tikzcd}
		s\mathbf{Set}_\bullet 
		\arrow[r] 
		\arrow[d, shift right=0.5ex, swap, "{f_\mathbb{Z}}"]& 
		\mathcal{H}_\bullet
		\arrow[r, shift left=0.5ex, "{\Sigma^\infty}"] 
		\arrow[d, shift right=0.5ex, swap, "{f_\mathbb{Z}}"] & 
		\mathcal{SH}
		\arrow[l, shift left=0.5ex, "{\Omega^\infty}"] \\
		C_+(\mathbf{Ab})
		\arrow[r] 
		\arrow[u, shift right=0.5ex, swap, "U"]  &
		\mathcal{H}_+(\mathbf{Ab})
		\arrow[u, shift right=0.5ex, swap, "U"] & 
	\end{tikzcd}
\end{equation}

Now, it is time to complete the stable part of this picture. As we saw earlier, for topological spaces, the suspension functor eventually became the shift functor. 
Also, the normalized chain complex of a simplicial abelian group after suspension is the normalized chain complex of the original spaces but shifted one time to the left. Therefore, it makes sense to use the shift functor, denoted by $[1]$, as the analogue of the suspension functor. When working with chain complexes of abelian groups with non-negative support, this functor is not an auto-equivalence. But, this problem can be addressed easily by expanding our category to include all complexes of abelian groups. Because in this category the shift functor has an obvious inverse, namely the shift-back functor $[-1]$. 

Let $A$ be an arbitrary complex of abelian groups. Define $A_{\geq -n}$ such that it is equal to $A$ for all $m\geq -n$ and $0$ elsewhere. By shifting this complex $n$ times, we get a complex $A_{\geq -n}[n]$ which is connective. We can recover $A$ from such complexes as
$$A = \text{colim}_n A_{\geq -n}[n][-n]$$
using the maps $\epsilon_n[-n-1]$ where $\epsilon_n$ is the inclusion $A_{\geq n}[n][1] \hookrightarrow A_{\geq n}[n+1]$. Therefore, we can construct arbitrary complexes of abelian groups using connective complexes and the inclusions $\epsilon_n$ for each integer $n$. 

\begin{remark}
	We saw that expanding our category to include all complexes of abelian groups is basically the same as considering sequences of objects in $C_+ (\mathbf{Ab})$ with the bonding maps $\epsilon_n[-n-1]$. Going back to the definition of spectra of spaces, we can see that the steps taken in both cases are exactly the same. This makes the analogy between the spaces and connective chain complexes of abelian groups more precise.
\end{remark}

The next piece of our puzzle is given by the fact that the free-forgetful adjunction raises all the way up to the spectra. And by taking derived functors, our picture of classical homotopy theory looks like the following: 

\begin{equation}
	\begin{tikzcd}
		s\mathbf{Set}_\bullet 
		\arrow[r] 
		\arrow[d, shift right=0.5ex, swap, "{f_\mathbb{Z}}"]& 
		\mathcal{H}_\bullet
		\arrow[r, shift left=0.5ex, "{\Sigma^\infty}"] 
		\arrow[d, shift right=0.5ex, swap, "{f_\mathbb{Z}}"] & 
		\mathcal{SH}
		\arrow[l, shift left=0.5ex, "{\Omega^\infty}"] 
		\arrow[d, shift right=0.5ex, swap, "{f_\mathbb{Z}}"] \\
		C_+(\mathbf{Ab})
		\arrow[r] 
		\arrow[u, shift right=0.5ex, swap, "U"]  &
		\mathcal{H}_+(\mathbf{Ab})
		\arrow[r, shift left=0.5ex, "{\Sigma^\infty}"] 
		\arrow[u, shift right=0.5ex, swap, "U"] & 
		D(\mathbf{Ab})
		\arrow[l, shift left=0.5ex, "{\Omega^\infty}"] 
		\arrow[u, shift right=0.5ex, swap, "U"]
	\end{tikzcd}
\end{equation}

Based on this picture, we can ask what the free-forgetful adjunction means between $\mathcal{SH}$ and $D(\mathbf{Ab})$. It turns out that this adjunction identifies $D(\mathbf{Ab})$ with the homotopy category of $H\mathbb{Z}\text{-modules}$ in $\mathcal{SH}$. In fact, we need to introduce the notion of symmetric spectra for this to be true - see \cite{Barns2020} for more information about symmetric spectra. The spectrum $H\mathbb{Z}$ lifts to a commutative ring object in the category of symmetric spectra which allows us to define $H\mathbb{Z}\text{-modules}$. An $H\mathbb{Z}\text{-module}$ is a symmetric spectrum $E$ with an action $H\mathbb{Z} \wedge E \to E$, satisfying the usual module conditions. So, the forgetful functor becomes the functor that forgets the $H\mathbb{Z}\text{-module}$ structure and the free functor becomes the free $H\mathbb{Z}\text{-module}$ functor. Therefore, we can change our picture to the following: 

\begin{equation}\label{classicalpic}
	\begin{tikzcd}
		\mathbf{Set}_\bullet 
		\arrow[hookrightarrow]{r} 
		\arrow[d, shift right=0.5ex, swap, "{f_\mathbb{Z}}"] &
		s\mathbf{Set}_\bullet 
		\arrow[r] 
		\arrow[d, shift right=0.5ex, swap, "{f_\mathbb{Z}}"]& 
		\mathcal{H}_\bullet
		\arrow[r, shift left=0.5ex, "{\Sigma^\infty}"] 
		\arrow[d, shift right=0.5ex, swap, "{f_\mathbb{Z}}"] & 
		\mathcal{SH}
		\arrow[l, shift left=0.5ex, "{\Omega^\infty}"] 
		\arrow[d, shift right=0.5ex, swap, "{H\mathbb{Z} \wedge -}"] \\
		\mathbf{Ab}
		\arrow[hookrightarrow]{r}
		\arrow[u, shift right=0.5ex, swap, "U"] &
		C_+(\mathbf{Ab})
		\arrow[r] 
		\arrow[u, shift right=0.5ex, swap, "U"]  &
		\mathcal{H}_+(\mathbf{Ab})
		\arrow[r, shift left=0.5ex, "{\Sigma^\infty}"] 
		\arrow[u, shift right=0.5ex, swap, "U"] & 
		D(\mathbf{Ab})
		\arrow[l, shift left=0.5ex, "{\Omega^\infty}"] 
		\arrow[u, shift right=0.5ex, swap, "H\mathbb{Z}"]
	\end{tikzcd}
\end{equation}
Here, we have added pointed sets and abelian groups to the picture for reasons that will become clear in later sections. Having this picture in mind, we will track the development of the Motivic Homotopy Theory in the next two sections.

\section{The picture in motivic homotopy theory: Part 1}\label{M1}

Motivic Homotopy Theory, in the most general case, wants to build a homotopy theory for the category $Sm/S$ of smooth schemes of finite type over a Noetherian scheme $S$. As this category does not have the properties needed to become a model category, we can not do homotopy theory in itself. For example, it does not contain all the small colimits. Therefore, we should find a bigger category that has the necessary properties and contains $Sm/S$. One way to address the problem of not having colimits is to embed $Sm/S$ into the category of presheaves on it $\text{PSh}(Sm/S)$ using the Yoneda embedding, but this approach is not good enough either. As mentioned in \cite{Voevodsky1998}, if we consider a covering $X = U \cup V$ of a scheme $X$ by two open subsets with respect to the Zariski topology and if we define our category of spaces to be $\text{PSh}(Sm/S)$, the union $U \cup_{\text{PSh}} V$ is not the same as $X$. To solve the problem of finding a nice category of spaces, Morel and Voevodsky used the category of Nisnevich simplicial sheaves on $Sm/S$. This category is equipped with a model structure given by André Joyal in \cite{Joyal1984}. This approach has a big positive point as it allows us to transfer the structures available for simplicial sets to our new category of spaces. Also, this choice is closely related to the universal homotopy category of $Sm/S$ - a notion introduced by Daniel Dugger in \cite{Dugger2001}. In his paper, Dugger proves that the homotopy theory built by Morel and Voevodsky is a result of performing two localizations to the universal model category of $Sm/S$. One of these localizations is with respect to the affine line $\mathbb{A}^1$, which plays the role of the unit interval $[0,1]$ in motivic homotopy theory. This localization would be addressed more precisely later on. The other localization comes from the Nisnevich topology on $Sm/S$. The choice of this topology is because while it is strictly stronger than the Zariski topology and strictly weaker than the étale topology, it inherits nice properties of both. To mention a few 
\begin{enumerate}
	\item 
	similar to the Zariski topology, the Nisnevich cohomological dimension of a scheme of Krull dimension $d$ is equal to $d$;
	\item
	similar to the étale topology, in the Nisnevich topology smooth pairs $(X, Y)$ are locally equivalent to pairs of the form $(\mathbb{A}^n, \mathbb{A}^m)$;
	\item 
	similar to the Zariski topology, algebraic $K\text{-theory}$ has Nisnevich descent;
	\item 
	similar to the étale topology, we can use Čech cochains to compute Nisnevich cohomology.
\end{enumerate}

In this section, we first introduce Joyal's model structure on a site and discuss how left Bousfield localizations are performed on this model category in Subsection \ref{41}. Then, we turn our focus to the site $(Sm/k)_\text{Nis}$, where $k$ is a field of characteristic zero, and perform a localization with respect to the affine line in the model structure from previous subsection. This results in the motivic unstable homotopy category. After this, we find our way to the motivic stable homotopy category and finish the first part of our motivic picture in \ref{43}. The two main references for this section are  \cite{Voevodsky1998} and \cite{MorelVoevodsky1999}.

\subsection{Homotopy category of sites with an interval}\label{41}

Let $T$ be an essentially small site with enough points. As mentioned earlier, we want to work with the category of simplicial sheaves on $T$, denoted by $\Delta^{\text{op}}\text{Shv}(T)$. The most important thing about this category is that it contains a version of $T$ and a version of simplicial sets in it. So, somehow, we are putting both of these categories together to create a platform for doing homotopy theory. The Yoneda embedding gives us a functor $T\to \text{PShv}(T)$. For the categories that we work with, representable presheaves are sheaves, so we get a fully faithful embedding $T\to \text{Shv}(T)$. Then, from a sheaf of sets, we can construct a simplicial object, and get a simplicial sheaf on $T$. So, every object $X\in T$ gives us a simplicial sheaf on $T$, and therefore $\Delta^{\text{op}}\text{Shv}(T)$ has a version of $T$ in it. On the other hand, every simplicial set defines a constant simplicial presheaf on $T$. By considering its associated sheaf, we can see that $\Delta^{\text{op}}\text{Shv}(T)$ has a version of $s\mathbf{Sets}$ in it. We call  $\Delta^{\text{op}}\text{Shv}(T)$ the category of spaces. The same steps above can be done for pointed simplicial sheaves on $T$. As a result, we get the category of pointed spaces, denoted by $\Delta^{\text{op}}\text{Shv}_\bullet(T)$. There is an evident free-forgetful adjunction between spaces and pointed spaces. 

Working with simplicial sheaves has a big advantage for us. Constructions for simplicial sets extend section-wise to our spaces, including limits, colimits, etc. Here, we discuss some of them. Let $(X, x)$ and $(Y, y)$ be two pointed spaces. The wedge sum of them is the sheaf associated with the presheaf 
$$ ((X, x) \vee (Y, y))(U):= (X, x_U) (U) \vee (Y, y_U)(U).$$
The smash product of them is the sheaf associated with the presheaf
$$((X, x) \wedge (Y, y))(U):= (X, x_U) (U) \wedge (Y, y_U)(U).$$
The internal hom objects in the category of pointed spaces are given by the right adjoint of the smash product. The suspension functor is defined as the smash product with the simplicial sheaf $S^1$, where $S^1$ is the simplicial set $\Delta^1 / \partial \Delta^1$. The right adjoint of the suspension functor gives the loop space functor.

In a letter to Alexander Grothendieck, \cite{Joyal1984}, André Joyal introduced a model structure for this category which uses point-wise weak equivalences. Before we define this model structure, we need to recall the definition of a point of $T$.

\begin{definition}
	A \textit{point} of $T$ is a functor $x:\text{Shv}(T)\to \mathbf{Sets}$, that commutes with all colimits and finite limits.
\end{definition}

\begin{definition}\label{we}
	A morphism of spaces $f: X\to Y$ is called a \textit{point-wise weak equivalence} if for any point $x$ of $T$ the morphism of simplicial sets $x^*(f): x^*(X)\to x^*(Y)$ is a weak equivalence.
\end{definition}

\begin{theorem}[Joyal \cite{Joyal1984}]\label{joyal}
	The category of spaces with the classes of 
	\begin{enumerate}
		\item 
		Weak equivalences: the point-wise weak equivalences from definition \ref{we}
		\item 
		Cofibrations: monomorphisms
	\end{enumerate}
	is a proper closed simplicial model category. 
\end{theorem}

Following Morel and Voevodsky, we call this model structure the simplicial model structure and denote its homotopy category with $\mathcal{H}_s(T)$. The canonical functor from simplicial sets to this model category preserves weak equivalences.

A similar theorem can be stated for pointed spaces. Hence we get the pointed homotopy category denoted by $\mathcal{H}_{s,\bullet}(T)$. The free-forgetful adjunction preserves weak equivalences in both directions, therefore we get an adjunction between $\mathcal{H}_s(T)$ and $\mathcal{H}_{s,\bullet}(T)$.


Now we have a model structure for our category of spaces. However, this model structure would not be satisfying in many cases, as usually there would be other morphisms in $T$ that we expect to be weak equivalences. So, we need a procedure that allows us to invert a class of morphisms and still get a model category. This procedure is called (left) Bousfield localization. Let $A$ be a set of morphisms that we want to invert. First, we need to notice that only inverting the morphisms in $A$ would not be enough. For example, for a morphism $f: X\to Y$ in $A$ we expect the induced map $X\times Z\to Y\times Z$ to be inverted as well. So, we need to expand $A$ to a set of morphisms we call $A\text{-local}$ morphisms.

\begin{definition}
	An object $X\in \mathcal{H}_s(T)$ is called $A\textit{-local}$ if for every $Y \in \mathcal{H}_s(T)$ and for every morphism $f:Z_1\to Z_2\in A$ the induced map
	$$\mathrm{Hom}_{ \mathcal{H}_s(T)}(Y\times Z_2, X)\to \mathrm{Hom}_{ \mathcal{H}_s(T)}(Y\times Z_1, X)$$
	is a bijection.
\end{definition}

\begin{definition}\label{Awe}
	A map of spaces $f:X\to Y$ is called an $A\textit{-weak}$ \textit{equivalence} if for every $A\text{-local}$ object $Z\in \mathcal{H}_s(T)$ the induced map 
	$$\mathrm{Hom}_{ \mathcal{H}_s(T)}(Y, Z)\to \mathrm{Hom}_{ \mathcal{H}_s(T)}(X, Z)$$
	is a bijection.
\end{definition}

Having defined $A\text{-weak}$ equivalence, (left) Bousfield localization gives us the following model structure.

\begin{theorem}[\cite{MorelVoevodsky1999}, Theorem 2.5]
	The category of spaces with the classes of 
	\begin{enumerate}
		\item 
		Weak equivalences: the $A\text{-weak}$ equivalences of definition \ref{Awe}
		\item 
		Cofibrations: monomorphisms
	\end{enumerate}
	is a model category, denoted by $\mathcal{H}(T, A)$. The inclusion functor $\mathcal{H}(T, A)\to \mathcal{H}_s(T)$ has a left adjoint which identifies $\mathcal{H}(T, A)$ with the localization of $\mathcal{H}_s(T)$ with respect to $A\text{-weak}$ equivalences.
\end{theorem}

We call this model structure the $A\text{-model}$ structure. A similar theorem can be stated for pointed spaces. Hence, we get the pointed homotopy category denoted by $\mathcal{H}_{\bullet}(T, A)$. Again the free-forgetful adjunction gives us an adjunction between these two homotopy categories.

In motivic homotopy, we want the affine line $\mathbb{A}^1$ to play the role of the unit interval $[0,1]$. This is because some cohomology theories like algebraic $K\text{-theory}$ (over a regular base) and motivic cohomology (over a field of characteristic zero) are $\mathbb{A}^1\text{-invariant}$, just like topological $K\text{-theory}$ and singular cohomology (and other reduced cohomology theories) are $[0,1]\text{-invariant}$. So, in the general case of the homotopy theory of a site $T$, we first need to define what an interval is and discuss the localization with respect to an interval.

\begin{definition}
	Let $T$ be a site and $*$ be the final object of $\text{Shv}(T)$. An \textit{interval} of $T$ is a sheaf $I$ with the morphisms 
	\begin{align*}
		\mu: I\times I\to I \\
		i_0, i_1:*\to I
	\end{align*}
	that satisfy the following conditions:
	\begin{enumerate}
		\item 
		For the canonical morphism $I\to *$, 
		\begin{align*}
			\mu(i_0\times I) = \mu(I\times i_0) = i_0 p \\
			\mu(i_1\times I) = \mu(I\times i_1) = \text{id}
		\end{align*}
		\item 
		The morphism $i_0\coprod i_1: *\coprod *\to I$ is a monomorphism.
	\end{enumerate}
\end{definition}

Localization with respect to $I$ is defined by considering the $A\text{-model}$ structure for $A:= \{i_0: *\to I\}$. We call this model structure the $I\text{-model}$ structure. This model structure is proper by Theorem 3.2 of \cite{MorelVoevodsky1999}.

\begin{example}
	Let $\Delta$ be the category that has the finite sets $\{0, 1, ..., n\}$ as its objects and order-preserving functions as it morphisms. Consider the trivial Grothendieck topology for this category. If we consider the interval $\Delta^1$ for this site, the resulting homotopy category would be the homotopy category of simplicial sets.
\end{example}

\subsection{Motivic unstable homotopy category}\label{42}

In this subsection, we will apply the theorems from section \ref{41} to the category of smooth schemes of finite type over $k$, where $k$ is a field of characteristic zero. This category is denoted by $Sm/k$. In order to do so, we need to choose a topology and an interval for this category. 

The chosen topology is the Nisnevich topology for reasons that were explained at the beginning of this Section. Now, let's define the Nisnevich topology.

\begin{definition}\label{421}
	A finite family of étale coverings $\{U_i\to X\}_{i\in I}$ is called a \textit{Nisnevich covering}, if and only if, for every $x\in X$, there exists an $i$ and $u\in U_i$ over $x$, such that the induced map on residue fields is an isomorphism. These coverings define a pretopology on $Sm/k$. The corresponding topology is called the \textit{Nisnevich topology}. We will denote this site by $(Sm/k)_{\mathrm{Nis}}$.
\end{definition}	

The chosen interval in $(Sm/k)_{\mathrm{Nis}}$ is the affine line $\mathbb{A}^1$ for reasons that were explained towards the end of Subsection \ref{41}. Let $(Sm/S, \mathbb{A}^1)_{\mathrm{Nis}}$ be the above site with the interval $\mathbb{A}^1$. Based on section \ref{41}, our category of spaces, which we call the category of motivic spaces and denote by $\textbf{Spc}(k)$, is the category of Nisnevich simplicial sheaves on $Sm/k$, and the homotopy category of motivic spaces, denoted by $\mathcal{H}(k)$ is the homotopy category corresponding to the $\mathbb{A}^1\text{-model}$ structure on $\textbf{Spc}(k)$. The same can be done for pointed spaces where we get $\mathbf{Spc}_\bullet(k)$ and $\mathcal{H}_\bullet(k)$. So, up to this point, our picture of motivic homotopy theory in contrast with the classical picture looks like the following: 

\begin{equation*}\label{locpic}
	\begin{tikzcd}
		\Delta^{\text{op}}\mathrm{Shv}_{\text{Nis}, \bullet}(Sm/k)
		\arrow[r] & 
		\mathcal{H}^{\mathbb{A}^1}_{\text{Nis}, \bullet}(k).
	\end{tikzcd}
\end{equation*}

\begin{remark}
In Picture \ref{locpic}, the functor is given by localizing with respect to the affine line. But based on the paper \cite{Dugger2001}, we can take $\Delta^{\text{op}}\mathrm{PSh}_{\bullet}(Sm/k)$ as our category of motivic spaces and then perform two localizations instead of one: one with respect to the affine line and one with respect to the homotopy-colimit-type relations coming from the Nisnevich topology. As a result, we might use either of these categories in our big picture..
\end{remark}

\begin{remark}
There are other Quillen equivalent model structures for motivic spaces, such as one due to Jardine in \cite{Jardine2000} and one in \cite{Dundas2003}. We will freely use any of these model categories when needed.
\end{remark}

\subsection{Spheres, homotopy sheaves, and suspensions}\label{43}

The next step for us is to build the motivic stable homotopy category. To do that, we need to define spectra, and therefore, we need to define the suspension functor. In the classical setting, suspension was defined using the smash product with the circle $S^1$. But, what is the analogue of $S^1$, and in general, what is a sphere in the motivic setting? As we know, $\mathbf{Spc}_\bullet(k)$ has both pointed simplicial sets and $Sm/k$ in it. This results in having two different circles; one coming from the topological world of simplicial sets and one from the algebraic world of $Sm/k$. The topological circle is the simplicial sheave associated with the pointed simplicial set $S^1$. The algebraic circle is $\mathbb{G}_m = \mathbb{A}^1- \{0\}$ which is a pointed space with respect to the base point $1$. This is called the Tate circle. So, now the question is which one of these circles should be considered as the right analogue of the sphere $S^1$. Before answering this question, let's find out about the possible analogues of other spheres too. In order to do so, we first work over the complex numbers. Let $X$ be a smooth scheme over complex numbers. Considering its complex points we get a topological space that we denote with $X(\mathbb{C})$. The functor $X\mapsto X(\mathbb{C})$ extends to a realization functor $\mathbf{Spc}_\bullet(\mathbb{C}) \to \mathbf{Top}_\bullet$, which preserves weak equivalences. Now, let's see what kinds of objects are sent to spheres by this functor. Here is a list of such spaces:

\begin{enumerate}
	\item 
	The simplicial sheave $S^n = (S^1)^{\wedge n}$ as a motivic space is sent to $S^n$, in particular $S^1$ is sent to $S^1$;
	\item 
	The space $(\mathbb{A}^n - \{0\}, (1,1, ..., 1))$ is sent to spheres, in particular $\mathbb{G}_m$ is sent to an sphere;
	\item 
	The projective line $(\mathbb{P}^1, \infty)$ is sent to a sphere;
	\item 
	The Thom space (which would be defined later on when stating the homotopy purity theorem), denoted by $\mathrm{Th}(\mathbb{A}^n)$, is sent to a sphere.
\end{enumerate}
Looking at these potential spheres would make things confusing at first. But, everything gets solved when we know that all of them are generated (up to $\mathbb{A}^1\text{-weak}$ equivalence) by the smash product of $S^1$ and $\mathbb{G}_m$. More precisely, 

\begin{align*}
	\mathbb{P}^1 \simeq S^1 \wedge \mathbb{G}_m \\
	\mathbb{A}^n - \{0\} \simeq S^{n-1} \wedge \mathbb{G}_m^{\wedge n} \\
	\mathrm{Th}(\mathbb{A}^n) \simeq S^n \wedge \mathbb{G}_m^{\wedge n}.
\end{align*}
where $\mathrm{Th}(\mathbb{A}^n)$ is the Thom space. These weak equivalences are not limited to when we work over complex numbers and we have them for other bases as well. As a result of this observation, we define mixed spheres.

\begin{definition}
	The \textit{mixed sphere} $S^{p,q}$ is defined as $S^{p-q} \wedge \mathbb{G}_m^q$.
\end{definition}


\begin{definition}
	Let $(X, x_0)$ be a pointed motivic space. The \textit{bigraded homotopy sheaf} $\pi_{a,b}^{\mathbb{A}^1}(X, x_0)$ is defined as the sheaf associated with the presheaf 
	$$U\mapsto \mathrm{Hom}_{\mathcal{H}_\bullet(k)}(S^{a,b} \wedge U_+, (X, x_0)).$$
\end{definition}

We define suspension with respect to a few important spheres, but a similar definition goes for any mixed sphere.

\begin{align*}
	\Sigma_{S^1} (X, x_0) := S^1 \wedge (X, x_0) \\
	\Sigma_{\mathbb{G}_m} (X, x_0) := \mathbb{G}_m \wedge (X, x_0) \\
	\Sigma_{\mathbb{P}^1} (X, x_0) := \mathbb{P}^1 \wedge (X, x_0) \\
	\Sigma_{T} (X, x_0) := \mathrm{Th}(\mathbb{A}^1) \wedge (X, x_0)
\end{align*}
The reason is that when constructing the category of motivic spectra we need to stabilize with respect to $\Sigma_{S^1}$ and $\Sigma_{\mathbb{G}_m}$ both. But because there is a canonical isomorphism $\Sigma_{T} \cong \Sigma_{S^1} \circ \Sigma_{\mathbb{G}_m}$, we only need to consider the category of $T\text{-spectra}$. But, the same steps can be taken by considering the category of $\mathbb{P}^1\text{-spectra}$.

\subsection{Motivic stable homotopy category}\label{44}

\begin{definition}
A $T\textit{-spectrum}$ is a sequence of pointed motivic spaces $(E_0, E_1, E_2, ...)$ with a family of bonding maps $\{\epsilon_n: \Sigma_{T}E_n \to E_{n+1}\}_{n=0}^{\infty}$.
\end{definition}

\begin{example}
To each pointed motivic space $X$, we can correspond a $T\textit{-spectrum}$ $\Sigma_{T}^\infty X = (X, \Sigma_{\mathbb{P}^1}X, \Sigma_{T}^2 X, ...)$ with identity bonding maps. In particular, the $T\textit{-spectrum}$ 
$$\mathbb{S}_k = (k_+, T, T \wedge T, ...)$$
with identity bonding maps is the analogue of the sphere spectrum in the classical setting.
\end{example}

\begin{definition}\label{P1maps}
A map of $T\text{-spectra}$ $f: E\to F$ is a family of maps $\{f_n: E_n\to F_n\}_{n=0}^{\infty}$ that respects the bonding maps of both spectra i.e. the diagram 
\begin{equation*}
	\begin{tikzcd}
		\Sigma_{T} E_n
		\arrow[r, "{\Sigma_{T} f_n}"] 
		\arrow[d, "{\epsilon_n}"] & 
		\Sigma_{T} F_n
		\arrow[d, "{{\epsilon}^{'}_{n}}"] \\
		E_{n+1}
		\arrow[r, "F_{n+1}"] &
		F_{n+1}
	\end{tikzcd}
\end{equation*}
commutes for each integer $n\geq 0$.
\end{definition}

We denote the category of $T\text{-spectra}$ with maps given by Definition \ref{P1maps} by $\mathbf{Spt}_{T}(k)$. In the next step, we need to turn this category into a model category. We do that following Jardine in \cite{Jardine2000}. First, note that the notions of weak equivalence, fibration, and cofibration can be defined term-wise for $T\text{-spectra}$.

\begin{definition}\label{termwise}
Let $f: E\to F$ be a map of $T\text{-spectra}$.
\begin{enumerate}
\item
The map $f$ is called a \textit{term-wise weak equivalence} if, for every integer $n\geq 0$, the map $f_n: E_n\to F_n$ is a weak equivalence in the unstable sense i.e. an $\mathbb{A}^1\text{-weak}$ equivalence.
\item
The map $f$ is called a \textit{term-wise cofibration} if, for every integer $n \geq 0$, the map $f_n: E_n\to F_n$ is a cofibration in the unstable sense i.e. a monomorphism.
\item
The map $f$ is called a \textit{term-wise fibration} if, for every integer $n \geq 0$, the map $f_n: E_n\to F_n$ is a fibration in the unstable sense i.e. a map having right lifting property with respect to monomorphisms that are also an $\mathbb{A}^1\text{-weak}$ equivalence.
\end{enumerate}
\end{definition}

\begin{theorem}[\cite{Jardine2000}, Lemma 1.2]\label{termmodel}
\begin{enumerate}
\item 
The category $\mathbf{Spt}_{T}(k)$ with the classes of 
\begin{enumerate}
		\item 
		Weak equivalences: the term-wise weak equivalences from definition \ref{termwise}
		\item 
		Fibrations: the term-wise fibrations from definition \ref{termwise}
\end{enumerate}
becomes a proper simplicial model category.
\item
The category $\mathbf{Spt}_{T}(k)$ with the classes of 
\begin{enumerate}
		\item 
		Weak equivalences: the term-wise weak equivalences from definition \ref{termwise}
		\item 
		Cofibrations: the term-wise cofibrations from definition \ref{termwise}
\end{enumerate}
becomes a proper simplicial model category.
\end{enumerate}
\end{theorem}

Starting with the first model structure in theorem \ref{termmodel}, we need to extend the class of weak equivalences. To do that, we should first define stable homotopy sheaves of $T\text{-spectra}$ using our definition of bigraded homotopy sheaves of pointed spaces, just like what we did in the classical setting.

\begin{definition}
Let $E$ be a $T\text{-spectrum}$. The colimit of 
\begin{equation*}
	\begin{tikzcd}
		\pi_{a+2n, b+n}^{\mathbb{A}^1} (E_n) 
		\arrow[r, "\Sigma_{T}"] & 
		\pi_{a+2n+2, b+n+1}^{\mathbb{A}^1} (\Sigma_{T} E_n) 
		\arrow[r, "\sigma_{n, *}"] & 
		\pi_{a+2n+2, b+n+1}^{\mathbb{A}^1} (E_{n+1}) 
		\arrow[r] & 
		...
	\end{tikzcd}
\end{equation*}
is the motivic stable homotopy sheaf $\pi_{a,b}^{\mathbb{A}^1}(E)$.
\end{definition}

\begin{definition}\label{mswe}
A map of $T\text{-spectra}$ $f: E\to F$ is called a \textit{motivic stable weak equivalence} if the induced map 
$$\pi_{a,b}^{\mathbb{A}^1}(E)\to \pi_{a,b}^{\mathbb{A}^1}(F)$$ 
is an isomorphism for each pair of integers $a,b \in \mathbb{Z}$.
\end{definition}

\begin{theorem}[\cite{Jardine2000}, Theorem 9.2]
The category $\mathbf{Spt}_{T}(k)$ with the classes of 
\begin{enumerate}
		\item 
		Weak equivalences: the motivic stable weak equivalences from definition \ref{mswe}
		\item 
		Cofibrations: the term-wise cofibrations from the first model structure in theorem \ref{termmodel}
\end{enumerate}
becomes a proper simplicial model category. We call the resulting homotopy category the motivic stable homotopy category and denote it by $\mathcal{SH}(k)$.
\end{theorem}

Up until this point, our picture of motivic homotopy theory looks like 

\begin{equation}\label{mpic1}
	\begin{tikzcd}
		Sm/k
		\arrow[r, "Y"] &
		\mathrm{PSh}_\bullet(Sm/k) 
		\arrow[r] &
		\Delta^{\text{op}}\mathrm{PSh}_{\bullet}(Sm/k)
		\arrow[r] & 
		\mathcal{H}^{\mathbb{A}^1}_{\text{Nis}, \bullet}(k)
		\arrow[r, shift left=0.5ex, "{\Sigma_{T}^\infty}"] & 
		\mathcal{SH}(k) 
		\arrow[l, shift left=0.5ex, "{\Omega_{T}^\infty}"] 
	\end{tikzcd}
\end{equation}
where we have added $Sm/k$ and its Yoneda embedding into $\mathrm{PSh}_\bullet(Sm/k)$ to the picture and $\Omega_{T}^\infty$ is the right adjoint of $\Sigma_{T}^\infty$.

\subsection{Cohomology theories}\label{45}
Just like spectra of topological spaces, each $T\text{-spectrum}$ defines one homology and one cohomology theory.

\begin{definition}
	Let $E$ and $F$ be two $T\text{-spectra}$. The $E\textit{-homology}$ of $F$ is defined to be
	$$E_{a,b}(F) := \text{Hom}_{\mathcal{SH}(k)}(\Sigma^{a,b} \mathbb{S}_k, E\wedge F),$$
	and the $E\textit{-cohomology}$ of $F$ is defined to be 
	$$E^{a,b}(F) := \text{Hom}_{\mathcal{SH}(k)}(F, \Sigma^{a,b} E).$$
	For a pointed motivic space $X$, its $E\text{-homology}$ and $E\text{-cohomology}$ are respectively $E_{a,b}(\Sigma_{T}^\infty X)$ and $E^{a,b}(\Sigma_{T}^\infty X)$.
\end{definition}

Now, an important question is if $T\text{-spectra}$ represent all the cohomology theories for schemes over $k$, just like spectra did in the classical setting. We will only discuss the representability of motivic cohomology, the analogue of singular cohomology, only because it plays an important role in our picture of motivic homotopy, just like singular cohomology did in the picture of classical homotopy.

In the classical setting, singular cohomology was represented by the Eilenberg-MacLane spectrum, which was built out of Eilenberg-MacLane spaces. Here, we will briefly discuss the process of constructing Eilenberg-MacLane spaces and the Eilenberg-MacLane spectrum in the classical setting, and then, we will take similar steps in the motivic setting to build the motivic Eilenberg-MacLane $T\text{-spectrum}$.

Let $X$ be a pointed topological space, and define the $n\text{-th}$ symmetric product of $X$, denoted by $\text{Sym}^n X$, to be the quotient $X^n / S_n$, where the symmetric group $S_n$ acts on $X^n$ by permuting the factors. Now, we can define the infinite symmetric product of $X$ as 
$$\text{Sym}^\infty X = \text{colim}_n \text{ Sym}^n X$$
where the maps are given by the inclusions $\text{Sym}^n X \hookrightarrow \text{Sym}^{n+1} X$. The Dold-Thom theorem gives a way of computing the homotopy groups of $\text{Sym}^\infty X$.

\begin{theorem}[Dold-Thom \cite{DoldThom1958}]
Let $(X, x)$ be a pointed topological space. Then, there is a natural isomorphism 
$$\pi_n(\mathrm{Sym}^\infty X, x) \cong H_n(X, x)$$
where $H$ is the reduced singular homology with coefficients in $\mathbb{Z}$.
\end{theorem}

Having this theorem, we know that $\text{Sym}^\infty S^n$ is the Eilenberg-MacLane space $K(\mathbb{Z}, n)$. Therefore, the Eilenberg-MacLane spectrum $H\mathbb{Z}$ is given by $(\text{Sym}^\infty S^0, \text{Sym}^\infty S^1, \text{Sym}^\infty S^2, ...)$.

To construct the motivic Eilenberg-MacLane we will substitute $S^1$ for $\mathbb{P}^1$, and more generally, we will replace $S^n$ with the mixed sphere $S^{2n,n} = (\mathbb{P}^1)^n$. So, the $n\text{-th}$ symmetric product of $\mathbb{P}^1$ is
$$\text{Sym}^\infty \mathbb{P}^1 := (\mathbb{P}^1)^n / S_n$$
and the infinite symmetric product of $\mathbb{P}^1$ is 
$$\text{Sym}^\infty \mathbb{P}^1 = \text{colim}_n \text{ Sym}^n \mathbb{P}^1$$
where the maps are given by inclusion. Now, we define the motivic Eilenberg-MacLane $T\text{-spectrum}$.

\begin{definition}\label{motcoh}
The motivic Eilenberg-MacLane $T\text{-spectrum}$, $M\mathbb{Z}$ is 
$$(\text{Sym}^\infty S^{0,0}, \text{Sym}^\infty S^{2,1}, \text{Sym}^\infty S^{4,2}, ...).$$
\end{definition}
	
One thing to notice here is that the motivic Eilenberg-Maclane $T\text{-spectrum}$ is not a perfect analogue of the Eilenberg-MacLane spectrum, as it does not share exactly the same properties. For example, the Eilenberg-MacLane space $K(G, n)$ has only one non-trivial homotopy group which is $\pi_n(K(G, n)) = G$, but the analogous result is not true for the motivic Eilenberg-MacLane space $\text{Sym}^\infty S^{2n,n}$. However, it represents motivic cohomology, which is the result we were after.
\begin{theorem}[Voevodsky \cite{Voevodsky1998}]
Over a field $k$ of characteristic zero, motivic cohomology is represented by 
\begin{align*}
H^{p,q} (X, \mathbb{Z}) = M\mathbb{Z}^{p,q}(\Sigma_{T}^\infty X) = \text{Hom}_{\mathcal{SH}(k)}(\Sigma_{T}^\infty X, \Sigma^{p,q} M\mathbb{Z})
\end{align*}
\end{theorem}

\section{The picture in motivic homotopy theory: Part 2}\label{M2}

Recall that the Eilenberg-Maclane spectrum provided us with a free-forgetful adjunction in the classical setting that led us to connect the stable homotopy category to the (unbounded) derived category of abelian groups. Using that adjunction, $D(\mathbf{Ab})$ was identified with the category of $H\mathbb{Z}\text{-modules}$ in the category of symmetric spectra. On the other hand, $M\mathbb{Z}$ lifts to a commutative ring object in the category of motivic symmetric $T\text{-spectra}$ which allows us to define $M\mathbb{Z}\text{-modules}$. Now, the question is to find out what category does this lead us to. In other words, what is the analogue of $D(\mathbf{Ab})$ in the motivic setting? It turns out that the analogue is Voevodsky's triangulated category of motives - see \cite{Voevodsky2000} for a full discussion of Voevodsky's motives. In this section, we will construct the second part of our motivic picture having the Diagram \ref{classicalpic} of classical homotopy in mind. The main reference for this section is \cite{RondigsOstvaerMain}.

\subsection{Motivic spaces with transfers}

First of all, in the second part of our motivic picture, we expect to work with abelian or additive categories, and as we know $Sm/k$ is neither. So, in the first step, we embed $Sm/k$ into the category of finite correspondences, denoted by $Cor_k$. This category has the same objects as $Sm/k$, and for morphisms it has 
\begin{equation*}
Cor_k(X, Y):= \text{free abelian group on elementary correspondences.} 
\end{equation*}
where an elementary correspondence is a subscheme $W$ of $X\times Y$ whose projection $W\to X$ is finite and onto a component of $X$. The category $Sm/k$ embeds into $Cor_k$ as the graph of each morphism $X\to Y \in Sm/k$ is an elementary correspondence. 

Now, the Yoneda embedding enables us to embed $Cor_k$ into the category of contravariant additive functors from $Cor_k$ to the category of abelian groups. This category is called the category of presheaves with transfers on $Cor_k$ and is denoted by $\text{PSh}^{\text{tr}}(Cor_k)$. 

Analogous to the embedding of abelian groups into simplicial abelian groups in the classical picture, we can embed presheaves with transfers into the category of $\Delta^{\text{op}}\text{PSh}^{\text{tr}}(Cor_k)$, which we call the category of motivic spaces with transfers. Notice that, we have an evident forgetful functor $U: \text{PSh}^{\text{tr}}(Cor_k) \to \mathrm{PSh}_\bullet(Sm/k)$ induced by the graph of $Sm/k \hookrightarrow Cor_k$. This functor has a left adjoint called the transfer functor denoted by $\mathbb{Z}^{\text{tr}}$. Moreover, we have a similar adjunction between $\Delta^{\text{op}}\text{PSh}^{\text{tr}}(Cor_k)$ and $\Delta^{\text{op}}\mathrm{PSh}_{\bullet}(Sm/k)$.

So up to this point, our motivic picture looks like the following:

\begin{equation}\label{mpic1}
	\begin{tikzcd}
		Sm/k
		\arrow[r, "Y"] 
		\arrow[d] &
		\mathrm{PSh}_\bullet(Sm/k) 
		\arrow[r] 
		\arrow[d, shift right=0.5ex, swap, "{\mathbb{Z}^{\text{tr}}}"] &
		\Delta^{\text{op}}\mathrm{PSh}_{\bullet}(Sm/k)
		\arrow[r] 
		\arrow[d, shift right=0.5ex, swap, "{\mathbb{Z}^{\text{tr}}}"] & 
		\mathcal{H}^{\mathbb{A}^1}_{\text{Nis}, \bullet}(k)
		\arrow[r, shift left=0.5ex, "{\Sigma_{\mathbb{P}^1}^\infty}"] & 
		\mathcal{SH}(k) 
		\arrow[l, shift left=0.5ex, "{\Omega_{\mathbb{P}^1}^\infty}"] \\
		Cor_k 
		\arrow[r, "Y"] &
		\text{PSh}^{\text{tr}}(Cor_k)
		\arrow[r] 
		\arrow[u, shift right=0.5ex, swap, "U"] &
		\Delta^{\text{op}}\text{PSh}^{\text{tr}}(Cor_k)
		\arrow[u, shift right=0.5ex, swap, "U"]
	\end{tikzcd}
\end{equation}

\subsection{Unstable homotopy category of motivic spaces with transfers}

Having defined motivic spaces with transfers, the following theorem leads us to its homotopy category. Here, we are using the model structure from \cite{Dundas2003} for motivic spaces, which is equivalent to the $\mathbb{A}^1\text{-model}$.

\begin{theorem}[Röndigs and Østær, \cite{RondigsOstvaerMain} Lemma 9 and Lemma 10]\label{MMtr}
There exists a monoidal and simplicial model structure for motivic spaces with transfers where weak equivalences (resp. fibrations) are weak equivalences (resp. fibrations) of underlying motivic spaces. Therefore, the forgetful functor $U$ detects and preserves motivic weak equivalences and fibrations.
\end{theorem}

In the classical picture, we identified the category of simplicial abelian groups with the category of chain complexes of abelian groups with non-negative support. In the motivic picture, too, we have an analogous Dold-Kan equivalence between $\Delta^{\text{op}}\text{PSh}^{\text{tr}}(Cor_k)$ and $C_+ (\text{PSh}^{\text{tr}}(Cor_k))$. Therefore, we may transport the model structure on $\Delta^{\text{op}}\text{PSh}^{\text{tr}}(Cor_k)$ to a model structure on $C_+ (\text{PSh}^{\text{tr}}(Cor_k))$ and get a Quillen equivalence. This enables us to extend to the (unbounded) derived category in later stages.

Based on the Theorem \ref{MMtr}, we can develop our picture even further and get the following diagram:

\begin{equation}\label{mpic1}
	\begin{tikzcd}
		Sm/k
		\arrow[r, "Y"] 
		\arrow[d] &
		\mathrm{PSh}_\bullet(Sm/k) 
		\arrow[r] 
		\arrow[d, shift right=0.5ex, swap, "{\mathbb{Z}^{\text{tr}}}"] &
		\Delta^{\text{op}}\mathrm{PSh}_{\bullet}(Sm/k)
		\arrow[r] 
		\arrow[d, shift right=0.5ex, swap, "{\mathbb{Z}^{\text{tr}}}"] & 
		\mathcal{H}^{\mathbb{A}^1}_{\text{Nis}, \bullet}(k)
		\arrow[r, shift left=0.5ex, "{\Sigma_{\mathbb{P}^1}^\infty}"] 
		\arrow[d, shift right=0.5ex, swap, "{\mathbb{Z}^{\text{tr}}}"] & 
		\mathcal{SH}(k) 
		\arrow[l, shift left=0.5ex, "{\Omega_{\mathbb{P}^1}^\infty}"] \\
		Cor_k 
		\arrow[r, "Y"] &
		\text{PSh}^{\text{tr}}(Cor_k)
		\arrow[r] 
		\arrow[u, shift right=0.5ex, swap, "U"] &
		C_+ (\text{PSh}^{\text{tr}}(Cor_k))
		\arrow[u, shift right=0.5ex, swap, "U"] 
		\arrow[r]  &
		\mathcal{H}_+ (\text{PSh}^{\text{tr}}(Cor_k)).
		\arrow[u, shift right=0.5ex, swap, "U"]
	\end{tikzcd}
\end{equation}

\begin{remark}
In picture \ref{mpic1}, we have used the Dold-Kan equivalence to substitute the category of motivic spaces with transfers and its homotopy category with $C_+ (\text{PSh}^{\text{tr}}(Cor_k))$ and its homotopy category which we denote by $\mathcal{H}_+ (\text{PSh}^{\text{tr}}(Cor_k))$. As this equivalence is a Quillen equivalence, we can use either in our picture, but to make the analogy between our classical and our motivic picture, we will use the latter.
\end{remark}

\subsection{Stable homotopy category of motivic spaces with transfers}

As you can already see, the steps taken to construct the second part of our motivic picture are very similar to the ones from the first part. For the stable part, too, we will take similar steps, but this time, something of importance will happen. The category occupying the bottom right place in our motivic picture will be equivalent to Voevodsky's triangulated category of motives, which is conjectured to be the derived category of mixed motives. 


In the previous subsection, we introduced the Quillen equivalent model categories $\Delta^{\text{op}}\text{PSh}^{\text{tr}}(Cor_k)$ and $C_+ (\text{PSh}^{\text{tr}}(Cor_k))$, and in this section, we will work with symmetric spectra of these two categories. As the steps are very similar to those seen before, we will not go into detail.

When constructing the motivic stable homotopy category, we used suspension with respect to $\mathbb{P}^1 \simeq T = \text{Th}(\mathbb{A}^1)$. Now, if we apply the translation functor $\mathbb{Z}^{\text{tr}}$ to $T$, we get a motivic space with transfers, and if we define symmetric spectra with respect to the suspension coming from this motivic space with transfers, we get the category of motivic symmetric spectra with transfers. We denote this category by $\textbf{MSS}^{\text{tr}}$. 

\begin{theorem}[\cite{RondigsOstvaerMain}, Theorem 11]\label{step1}
The category of motivic symmetric spectra with transfers acquires a monoidal stable model structure such that the forgetful functor detects and preserves weak equivalences and fibrations between stably fibrant objects.
\end{theorem}

We can take similar steps using $\mathbb{P}^1$ instead of $T$, getting the category $\textbf{MSS}^{\text{tr}}_{\mathbb{P}^1}$. Based on \cite{RondigsOstvaerTR}, the resulting model categories are Quillen equivalents through a zig-zag of strict symmetric monoidal Quillen equivalences.

\begin{proposition}[\cite{RondigsOstvaerMain}, Proposition 31]
There is a zig-zag of strict monoidal Quillen equivalences between $\mathbf{MSS}^{\mathrm{tr}}$ and $\mathbf{MSS}^{\mathrm{tr}}_{\mathbb{P}^1}$.
\end{proposition}

Using the Dold-Kan equivalence, we can take similar steps for $C_+ (\text{PSh}^{\text{tr}}(Cor_k))$. We can take symmetric spectra of objects in this category with respect to $\mathbb{Z}^{\text{tr}}(\mathbb{G}_m, 1)$. We denote the resulting category by $C_+\textbf{SS}_{\mathbb{G}_m}$. Based on Lemma 27 in \cite{RondigsOstvaerMain}, this category acquires a model structure with level-wise scheme-wise quasi-isomorphisms as weak equivalences and monomorphisms as cofibrations. Moreover, the shift functor is a Quillen equivalence, so we could work with the shifted $\mathbb{Z}^{\text{tr}}(\mathbb{G}_m, 1)[1]$ and $C_+\textbf{SS}_{\mathbb{G}_m [1]}$ too. We can also define $C_+ \textbf{SS}_{\mathbb{P}^1}$ in a similar fashion by suspending with respect to $\mathbb{Z}^{\text{tr}}(\mathbb{P}^1, 1)$. There exists a similar zig-zag of strict symmetric monoidal Quillen equivalences for these categories too.

\begin{proposition}[\cite{RondigsOstvaerMain}, Proposition 32]
There is a zig-zag of strict monoidal Quillen equivalences between $C_+ \mathbf{SS}_{\mathbb{G}_m [1]}$ and $C_+ \mathbf{SS}_{\mathbb{P}^1}$.
\end{proposition}

Finally, we can connect the last two propositions through the next one.

\begin{proposition}[\cite{RondigsOstvaerMain}, Theorem 35]
The Dold-Kan equivalence induces a lax symmetric monoidal (right) Quillen equivalence between $\mathbf{MSS}^{\mathrm{tr}}_{\mathbb{P}^1}$ and $C_+\mathbf{SS}_{\mathbb{P}^1}$.
\end{proposition}

\subsection{The link to Voevodsky's triangulated category of motives}

Consider the motivic cohomology $T\text{-spectrum}$ $M\mathbb{Z}$ which was introduced in Definition \ref{motcoh}. As mentioned in Section \S\ref{M1}, we expect motivic cohomology to be the analogue of singular cohomology in the motivic homotopy theory. In our picture of classical homotopy theory, singular cohomology was key to connecting the stable homotopy category of symmetric spectra to the (unbounded) derived category of abelian groups. It turns out that motivic cohomology provides us with a similar connection in the motivic setting. 

The $T\text{-spectrum}$ $M\mathbb{Z}$ lifts to a commutative ring object in the category of symmetric $T\text{-spectra}$, so we can define the notion of an $M\mathbb{Z}\text{-module}$ in this category. On the other hand, symmetric $T\text{-spectra}$ satisfy the monoid axiom in \cite{SchwedeShipley2000}, therefore, we have a model structure for the category of $M\mathbb{Z}\text{-modules}$ where weak equivalences (resp. fibrations) are weak equivalences (resp. fibrations) of underlying motivic symmetric $T\text{-spectra}$. Moreover, this model structure is proper and monoidal as stated in Proposition 38 in \cite{RondigsOstvaerMain}. 

\begin{theorem}[\cite{RondigsOstvaerMain}, Theorem 58]\label{step2}
Suppose $k$ is a field of characteristic zero. Then there is a strict symmetric monoidal Quillen equivalence between the category of $M\mathbb{Z}\text{-modules}$ and $\mathbf{MSS}^{\mathrm{tr}}$.
\end{theorem}

This theorem completes our picture of motivic homotopy theory, and in the last stage of this picture, the forgetful functor $U$ becomes the functor that forgets the $M\mathbb{Z}\text{-module}$ structure and $\mathbb{Z}^{\text{tr}}$ becomes its left adjoint which is $M\mathbb{Z} \wedge [-]$. Before drawing our completed picture, we need to mention one more result.

\begin{theorem}[\cite{RondigsOstvaerMain}, Theorem 1]
If $k$ is a field of characteristic zero, the homotopy category of $M\mathbb{Z}\text{-modules}$ is equivalent to Voevodsky’s big category of motives $D\mathcal{M}(k)$ of $k$. This equivalence preserves the monoidal and triangulated structures.
\end{theorem} 

\begin{proof}
Combining Theorem \ref{step1} and Theorem \ref{step2} implies the result.
\end{proof}

As a result, our final picture of motivic homotopy theory looks like the following: 

\begin{equation}\label{mpic}
	\begin{tikzcd}
		Sm/k
		\arrow[r, "Y"] 
		\arrow[d] &
		\mathrm{PSh}_\bullet(Sm/k) 
		\arrow[r] 
		\arrow[d, shift right=0.5ex, swap, "{\mathbb{Z}^{\text{tr}}}"] &
		\Delta^{\text{op}}\mathrm{PSh}_{\bullet}(Sm/k)
		\arrow[r] 
		\arrow[d, shift right=0.5ex, swap, "{\mathbb{Z}^{\text{tr}}}"] & 
		\mathcal{H}^{\mathbb{A}^1}_{\text{Nis}, \bullet}(k)
		\arrow[r, shift left=0.5ex, "{\Sigma_{\mathbb{P}^1}^\infty}"] 
		\arrow[d, shift right=0.5ex, swap, "{\mathbb{Z}^{\text{tr}}}"] & 
		\mathcal{SH}(k) 
		\arrow[l, shift left=0.5ex, "{\Omega_{\mathbb{P}^1}^\infty}"]
		\arrow[d, shift right=0.5ex, swap, "{M\mathbb{Z} \wedge [-]}"] \\
		Cor_k 
		\arrow[r, "Y"] &
		\text{PSh}^{\text{tr}}(Cor_k)
		\arrow[r] 
		\arrow[u, shift right=0.5ex, swap, "U"] &
		C_+ (\text{PSh}^{\text{tr}}(Cor_k))
		\arrow[u, shift right=0.5ex, swap, "U"] 
		\arrow[r]  &
		\mathcal{H}_+ (\text{PSh}^{\text{tr}}(Cor_k))
		\arrow[u, shift right=0.5ex, swap, "U"] 
		\arrow[r, shift left=0.5ex, "{\Sigma^\infty}"] &
		D\mathcal{M}(k)
		\arrow[l, shift left=0.5ex, "{\Omega^\infty}"]
		\arrow[u, shift right=0.5ex, swap, "{M\mathbb{Z}}"]
	\end{tikzcd}
\end{equation}

\section{An introduction to the conjectural motivic $t\text{-structure}$}\label{T}

Up until now, we have used the word \lq\lq motivic" a lot, so a question that naturally arises is about the connection between our previous discussions and the conjectural category of mixed motives. In this section, we will introduce a conjectural bridge between Voevodsky's category of motives and the category of mixed motives. But before doing that, we need to introduce what a $t\text{-structure}$ is in a triangulated category.

\subsection{Definition of $t\text{-structure}$ and a fundamental example}

A simple instant of a $t\text{-structure}$ can be seen in the derived category of an abelian category. Let $\mathcal{A}$ be an abelian category. In its (unbounded) derived category $D(\mathcal{A})$, which is a triangulated category, we have a full subcategory $D(\mathcal{A})_{\geq 0}$, consisting of chain complexes with non-negative support. We can shift this full subcategory $-n$ times to get $D(\mathcal{A})_{\geq n}$, which consists of chain complexes with support $\geq n$. We can also define the full subcategory $D(\mathcal{A})_{\leq 0}$, consisting of chain complexes with non-positive support in a similar fashion, and shift it $-n$ times to get $D(\mathcal{A})_{\leq n}$. Having defined these subcategories, we have 
\begin{align*}
... \subseteq D(\mathcal{A})_{\leq -2}  \subseteq D(\mathcal{A})_{\leq -1} \subseteq D(\mathcal{A})_{\leq 0} \subseteq D(\mathcal{A})_{\leq 1} \subseteq D(\mathbf{Ab})_{\leq 2} \subseteq ... \\
... \subseteq D(\mathcal{A})_{\geq 2}  \subseteq D(\mathcal{A})_{\geq 1} \subseteq D(\mathcal{A})_{\geq 0} \subseteq D(\mathcal{A})_{\geq -1} \subseteq D(\mathbf{Ab})_{\geq -2} \subseteq ...,
\end{align*}
and for each integer $n$, each $A\in  D(\mathcal{A})_{\leq n}$, and each $B\in D(\mathcal{A})_{\geq n+1}$ we have
$$\text{Hom}_{D(\mathcal{A})}(A, B) = 0 \text{ and }  \text{Hom}_{D(\mathcal{A})}(B, A) = 0.$$
Also, for each $E\in D(\mathcal{A})$, there is a distinguished triangle 
$$E_0\to E\to E_1\to E_0[1]$$
where $E_0\in D(\mathcal{A})_{\leq 0}$ and $E_1\in D(\mathcal{A})_{\geq 1}$. Such a triangle also exists for each pair of integers $(n, n+1)$ instead of $(0,1)$. Lastly, the intersection 
$$D(\mathcal{A})_{\leq 0} \cap D(\mathcal{A})_{\geq 0}$$ 
which again is a full subcategory, is just the category $\mathcal{A}$ itself.

Generalizing the above data for arbitrary triangulated categories gives us the definition of a $t\text{-structure}$.

\begin{definition}
Let $\mathcal{C}$ be a triangulated category. A $t\textit{-structure}$ on $\mathcal{C}$ consists of two full subcategories $\mathcal{C}_{\geq 0}$ and $\mathcal{C}_{\leq 0}$ satisfying the following conditions:
\begin{enumerate}
\item
For $\mathcal{C}_{\geq n} := \mathcal{C}_{\geq 0}[n]$ and $\mathcal{C}_{\leq n} := \mathcal{C}_{\leq 0}[-n]$, we have $\mathcal{C}_{\leq -1} \subseteq \mathcal{C}_{\leq 0}$ and $\mathcal{C}_{\geq 1} \subseteq \mathcal{C}_{\geq 0}$.
\item
For any $X\in \mathcal{C}_{\leq 0}$ and any $Y\in \mathcal{C}_{\geq 1}$, we have $\mathrm{Hom}_{\mathcal{C}}(X, Y) = 0$.
\item
For any $Z\in \mathcal{C}$, there is a distinguished triangle
$$Z_0\to Z\to Z_1\to Z_0[1]$$
where $Z_0\in \mathcal{C}_{\leq 0}$ and $Z_1\in \mathcal{C}_{\geq 1}$.
\end{enumerate}
The \textit{heart} of this $t\text{-structure}$ is the full subcategory $\mathcal{C}_{\leq 0} \cap \mathcal{C}_{\geq 0}$.
\end{definition}

As in the fundamental example mentioned earlier, the heart of the canonical $t\text{-structure}$ of $D(\mathcal{A})$ gives us the abelian category $\mathcal{A}$ itself. In other words, a $t\text{-structure}$ provides us with a way to go back from $D(\mathcal{A})$ to $\mathcal{A}$. But an important thing to notice is that the could be many different $t\text{-structures}$ available on a triangulated category with different hearts, so working backwards is not always unique! 

The connection between $D\mathcal{M}(k)$ and the category of mixed motives is based on the existence of a conjectural $t\text{-structure}$ called the motivic $t\text{-structure}$.

\subsection{The conjectural motivic $t\text{-structure}$}

Voevodsky's triangulated category of motives satisfies properties expected from the derived category of mixed motives. So, a natural guess is that there should exist a $t\text{-structure}$, for which the heart is the category of mixed motives. Now, let's see a more precise statement of this conjecture based on \cite{Ayoub2017} by Joseph Ayoub. 

\begin{conjecture}
Let $D\mathcal{M}(k, A)$ be Voevodsky's category of motives with coefficients in $A$, and let $D(A)$ be the derived category of $A\text{-modules}$. Given an embedding $\iota: k\hookrightarrow \mathbb{C}$, let $\mathcal{T}_{\geq 0}^{\mathcal{M}}$ be the full subcategory of $D\mathcal{M}(k, A)$ consisting of those motives for which the Betti realization lands in $D(A)_{\geq 0}$. Recall that the Betti realization functor 
$$B_\iota: D\mathcal{M}(k, A)\to D(A)$$ 
is the realization functor that takes the motive of $X$ to the singular chain complex of $X^{\text{an}}$ with coefficients in $A$, where $X^{\text{an}}$ is the complex analytic space associated with $X$. Define $\mathcal{T}_{< 0}^{\mathcal{M}}$ to be the full subcategory whose objects are $N\in D\mathcal{M}(k, A)$ such that 
$$\text{Hom}_{D\mathcal{M}(k, A)}(M, N) = 0$$
for every $M\in \mathcal{T}_{\geq 0}^{\mathcal{M}}$. Then, the following properties should hold:
\begin{enumerate}
\item
The pair $(\mathcal{T}_{\geq 0}^{\mathcal{M}}, \mathcal{T}_{< 0}^{\mathcal{M}})$ defines a $t\text{-structure}$ on $D\mathcal{M}(k, A)$, which is independent of the choice of the complex embedding $\iota$.
\item
The Betti realization takes motives in $\mathcal{T}_{< 0}^{\mathcal{M}}$ to complexes in $D(A)_{< 0}$.
\item
Assuming that $A$ is a regular ring, this $t\text{-structure}$ can be restricted to the subcategory of geometric motives denoted by $D\mathcal{M}_gm (k, A)$.
\end{enumerate}
This $t\text{-structure}$, if exists, is called the motivic $t\text{-structure}$.
\end{conjecture}

Beilinson in \cite{Beilinson2010} proved that over a field of characteristic zero, the existence of the motivic $t\text{-structure}$ implies the standard conjectures on algebraic cycles. This result completes the connection between our motivic setting and the world of motives. Although, finding the motivic $t\text{-structure}$ in the most general case seems quite far, in the particular case of Tate motives over number fields it has been proved to exist by Marc Levine in \cite{Levine1993}.

\end{document}